\documentclass[12pt]{article}

\usepackage{graphicx,psfrag,color,amsmath}

\def\-{~-~}
\def\+{~+~}
\def\={~=~}
\def\eq{~\equiv~}
\def\beq{\begin{equation}}
\def\eeq{\end{equation}}
\def\beqn{\begin{eqnarray}}
\def\eeqn{\end{eqnarray}}
\def\bb{\begin{eqnarray*}}
\def\ee{\end{eqnarray*}}
\newcommand{\calle}[1]{(\ref{#1})}
\newcommand{\lp}{\left(}
\newcommand{\lb}{\left\lbrack}

\newcommand{\rp}{\right)}
\newcommand{\rb}{\right\rbrack}

\newcommand{\oma}[1]{\Omega_{#1}}
\newcommand{\omb}[1]{\Omega_{f#1}}
\newcommand{\omc}[1]{\Omega_{b#1}}
\newcommand{\kth}[1]{${#1}^{th}$}

\title{DOMINO WAVES}

\author{C.J. Efthimiou\thanks{Department of Physics,
                       University of Central Florida,
                       Orlando, FL 32816 ({\tt costas@physics.ucf.edu}).}
        \and
        M.D. Johnson\thanks{Department of Physics,
                       University of Central Florida,
                       Orlando, FL 32816.}
       }

\begin{document}

\maketitle

\begin{abstract}
Motivated by a proposal of Daykin \cite{Daykin}, we study the wave
that propagates along an infinite chain of dominoes and find the
limiting speed of the wave in an extreme case.\\[2mm]
\textbf{DOI}: 10.1137/S0036144504414505
\end{abstract}

 \pagestyle{myheadings}
 \thispagestyle{plain}
 \markboth{C.J. EFTHIMIOU AND M.D. JOHNSON}{DOMINO WAVES}

\section{Introduction}

Everyone is familiar with dominoes and has used them for fun.
A common game is to arrange the dominoes in a row and give a push
to the first. This generates a pleasing wave of falling dominoes.
The propagation happens at some speed $v$ (not necessarily constant).
A qualitative discussion  for general audiences is given by Walker
in \cite{Walker}.
Given the game's simplicity it is perhaps surprising to discover
that an exact
computation of the speed $v$ is quite difficult.
Daykin realized this in the following 1971 proposal \cite{Daykin}
to the readers of the \textsf{SIAM Review}:

\begin{center}
\footnotesize
\begin{minipage}{14cm}
  ``How fast do dominoes fall?" \\
   The ``domino theory of Southeast Asia" says that if Vietnam falls,
   then Laos falls,
   then Cambodia falls, and so on. Hearing a discussion of the theory
   led me to wonder
   about the proposed physical problem. The reader is invited to set
   his or her own
   ``reasonable" simplifying assumptions, such as perfectly elastic
   dominoes, constant
   coefficient of friction between dominoes and the table, and
   initial configuration
   with the dominoes equally spaced in a straight line, and so on.
\end{minipage}
\normalsize
\end{center}

In 1983, McLachlan et al. \cite{MLBCG} found a scaling law for the
speed $v$ in the limiting case of dominoes with zero thickness equally
spaced in a straight line.
With these assumptions the authors found the functional relation
 \beq
    v \= \sqrt{g\ell} \, G \! \lp {d\over \ell} \rp.
 \label{eq:1}
 \eeq
 Here $\ell$ is the height of the dominoes, $d$
the spacing between dominoes, and $G(x)$ an undetermined function
of $x$. This relation followed from dimensional analysis of the
problem\footnote{This paper should allow the reader to write down
a complete list  of assumptions needed to reach this conclusion.}.
McLachlan et al. proceeded to test the formula experimentally
using dominoes of heights $h=4.445cm$ and $h=8.890cm$. More
recently, Banks presented a simplified description of the effect
\cite{Banks}. His analysis, among other assumptions, assumes a
uniform propagation speed and conservation of momentum. A
different direction was taken by Shaw in a short paper \cite{Shaw}
describing how to model the domino effect as a computer simulation
and use it as an experiment in the undergraduate physics lab.

Equation \calle{eq:1} is not an entirely satisfactory solution, of
course, containing as it does an unknown function. This
recognition became the motivation for the present article. In this
work we develop an expression for the speed $v$. The result can be
cast as a particular forms of the scaling function $G$ arising
from a particular set of assumptions.

In order to highlight the basic physics behind the problem, we replace
the dominoes by massless rods topped with point masses $m$, as
seen in figure \ref{fig:pin-chain}.
A similar analysis in the case of dominoes with the shape of parallelepipeds
is straightforward, although there are some minor differences.

\begin{figure}[htb!]
\begin{center}
\psfrag{h}{$\ell$}
\psfrag{t}{$t$}
\psfrag{d}{$d$}
\includegraphics[width=15cm]{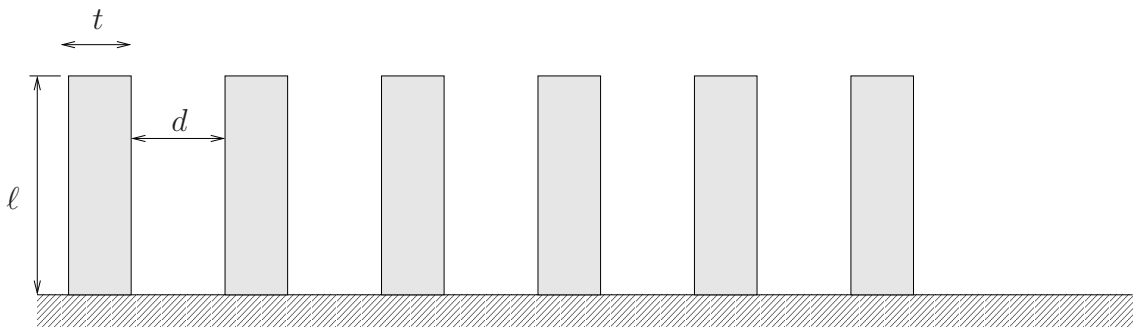}
\caption{A uniform chain of dominoes.}
\label{fig:uniform-chain}
\end{center}
\end{figure}

\begin{figure}[htb!]
\begin{center}
\psfrag{h}{$\ell$}
\psfrag{d}{$d$}
\psfrag{A1}{$A_1$}
\psfrag{A2}{$A_2$}
\psfrag{A3}{$A_3$}
\psfrag{A4}{$A_4$}
\psfrag{An}{$A_n$}
\includegraphics[width=12cm]{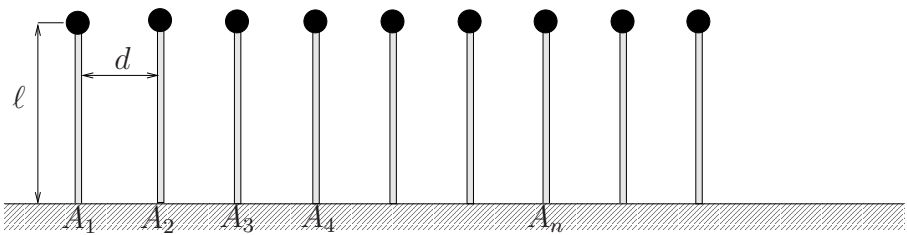}
\caption{A chain of massless rods carrying masses $m$ on top. The masses
         are indicated as finite spheres only for the sake of visualization.}
\label{fig:pin-chain}
\end{center}
\end{figure}

\section{The Model}

\subsection{The Assumptions}

We shall assume that:
\begin{enumerate}
\item \textit{The chain of  rods is uniform}. This means that all rods are identical and
are equally spaced along a straight line. Let $\ell$ be the length of a rod and
$m$ the mass on top.
\item \textit{The collisions are head-on}. This means that, seen from above,
      all the rods are and remain aligned on the same line. If this is not
      the case, then additional parameters are necessary to describe the
      collisions from point to point. The higher the asymmetry, the more
      parameters are needed, and the problem becomes highly complicated.
\item There is \textit{enough static friction} between the rods and the
      floor to keep the rods from sliding relative to the floor.
      Thus the rods pivot about
      fixed axes.
\item \textit{No energy is dissipated at the contact point}
      between the rods and the floor. This condition is independent of the
      previous one; it is possible for an object to rotate about an axis
      and yet to dissipate energy.
\item \textit{Collisions are instantaneous}. This means that the time
      interval $\Delta t$ during which a collision occurs is zero.
      The change in total angular momentum is then
      $\Delta\vec L=\int_0^{\Delta t}\vec \tau dt=0$.
      That is, total angular momentum
      is conserved during the collision, even though gravity is an external
      force and produces a non-zero torque.
\item \textit{The collisions are elastic}. This means that energy
is not
      dissipated during the collisions.
\item \textit{The rods are stiff}. This means that there is no deformation of
      the rods and thus no energy is converted to elastic potential energy of deformation.
      This condition is independent of the previous one; it is possible for
      the rods to be stiff and still dissipate energy during a collision.
\end{enumerate}

\subsection{Definition of Symbols}

To facilitate the calculations in the next section we present our notation
in advance. We label the rods sequentially with the numbers
$1,2,3,\dots,k,\dots,$ starting with 1 from the left end. We also label
$A_1,A_2,A_3,\dots,A_k,\dots,$ the pivot points of the rods as seen in figure
\ref{fig:pin-chain}. Then
\begin{itemize}
 \item $\theta_k$ is the angular displacement
       of a rod from the vertical, and $\omega_k=d\theta_k/dt$ is the corresponding
       time-dependent angular velocity.
\item $\oma{1}$ is the initial angular velocity of the first rod immediately after it is pushed.
 \item $\oma{k}$ is the initial angular velocity of the \kth{k} rod just as it begins to move.
This is of course the result of the collision with the \kth{(k-1)} rod (for $k>1$).
 \item $\omb{k}$ is the  angular velocity of the \kth{k} rod just before
       collision with the \kth{(k+1)} rod. (The subscript $f$ means `final'
       or `fallen'.)
 \item $\omc{k}$ is the  angular velocity of the \kth{k} rod just after
       collision with the \kth{(k+1)} rod
       ($b$ because the rod has just `bounced').
\item $\beta_1$ is the angle a rod forms with the vertical at the point of collision:
       $$
            \beta_1\=\sin^{-1}{d\over\ell}~.
       $$
  \item $T_k$ is the time the \kth{k} rod takes to fall from the vertical to its collision
  with the \kth{(k+1)} rod.
\end{itemize}
  In this notation, $\oma{k}$ is the angular velocity of the \kth{k} rod at
  $\theta_k=0$ and $\omb{k}$, $\omc{k}$ its angular velocities after it has fallen
  to $\theta=\beta_1$, just before and just after collision with the next rod
  respectively.

\begin{figure}[htbp!]
\begin{center}
\psfrag{a}{(a)}
\psfrag{B}{(b)}
\psfrag{c}{(c)}
\psfrag{d}{(d)}
\psfrag{e}{(e)}
\psfrag{f}{(f)}
\psfrag{b}{$\beta_1$}
\psfrag{h}{$\ell$}
\psfrag{t}{$t$}
\psfrag{A2}{$A_{k-1}$}
\psfrag{A3}{$A_k$}
\psfrag{A4}{$A_{k+1}$}
\psfrag{Oi4}{$\oma{k+1}$}
\psfrag{Oi3}{$\oma{k}$}
\psfrag{Of2}{$\omb{,k-1}$}
\psfrag{Of3}{$\omb{k}$}
\psfrag{O2}{$\omc{,k-1}$}
\psfrag{O3}{$\omc{k}$}
\includegraphics[height=15cm]{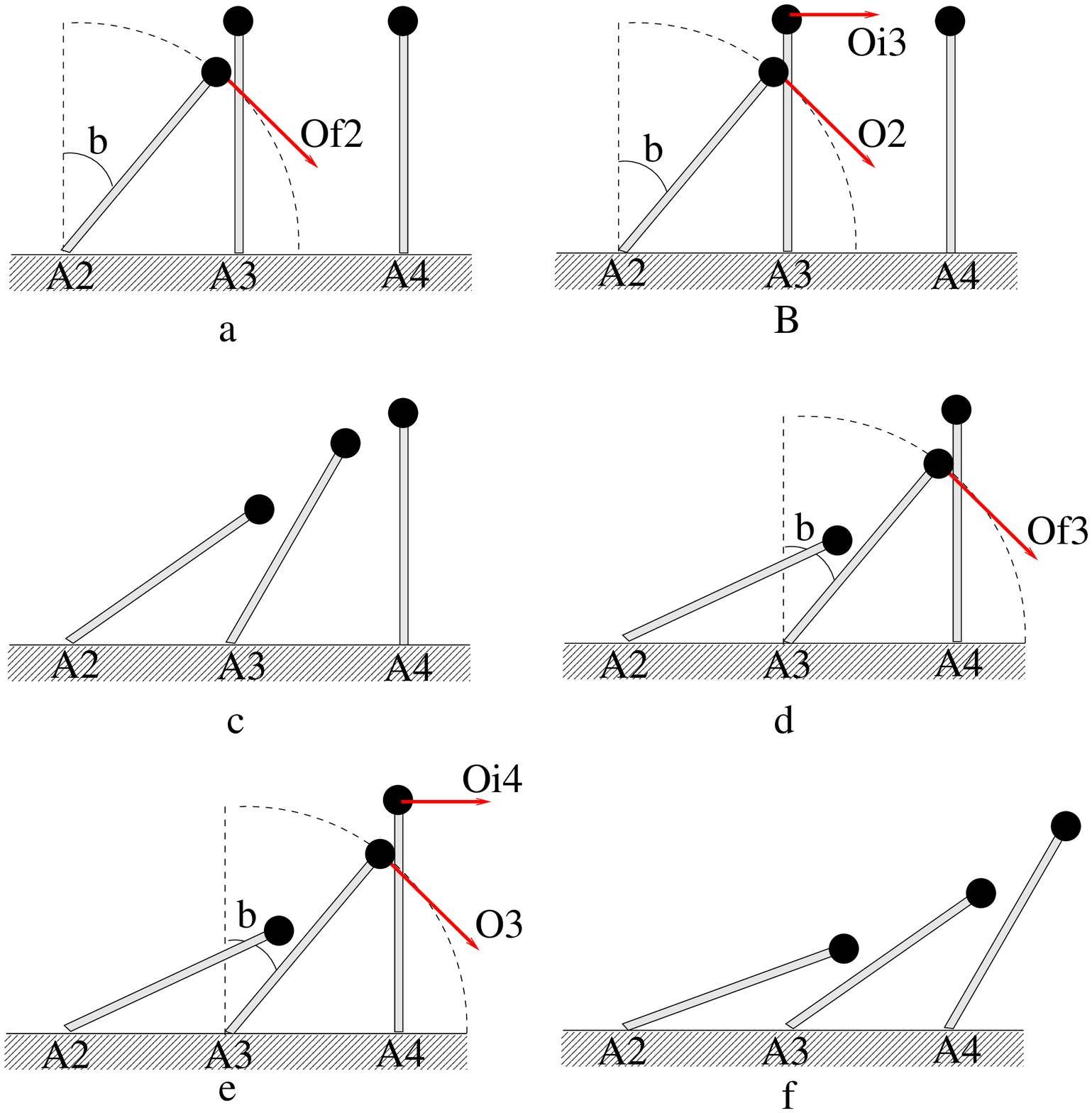}
\caption{The series of collisions between the \kth{(k-1)}, \kth{k}, and
         \kth{(k+1)} rods. In particular, the figure shows the state of the
         rods  \newline
        (a) just before the \kth{(k-1)} and \kth{k} rods collide.\newline
        (b) just after the \kth{(k-1)}  and \kth{k} rods  have collided.\newline
        (c) while the \kth{k} rod rotates towards the \kth{(k+1)} rod.\newline
        (d) just before the \kth{k} and \kth{(k+1)} rods collide.\newline
        (e) just after the \kth{k}  and \kth{(k+1)} rods  have collided.\newline
        (f) while the \kth{(k+1)} rod rotates towards the \kth{(k+2)} rod.
        }
\label{fig:domino}
\end{center}
\end{figure}

\subsection{Study of the Two-Rod Collisions}

Now examine the collision between the \kth{k} and \kth{(k+1)} rods. Our
assumptions guarantee that during the collision kinetic energy and
angular momentum are conserved, and that while a rod falls its total energy
(kinetic plus potential) is conserved. These conservation laws determine the solution.

Just before the collision the \kth{k} rod has angular velocity $\omb{k}$
and the \kth{(k+1)} rod is at rest. After collision the \kth{k} rod has angular
velocity $\omc{k}$ and the \kth{(k+1)} rod has angular
velocity $\oma{k+1}$.
Applying conservation of kinetic energy,
$$
  {1\over2}I\omb{k}^2\={1\over2}I\omc{k}^{2}
                             +{1\over2}I\oma{k+1}^2
$$
where $I=m\ell^2$ is the moment of inertia of a rod.
Therefore
\beq
  \omb{k}^2\=\omc{k}^{2}+\oma{k+1}^2~.
\label{eq:7}
\eeq

\begin{figure}[htb!]
\begin{center}
\psfrag{before}{BEFORE}
\psfrag{after}{AFTER}
\psfrag{h}{$\ell\,\cos\beta_1$}
\psfrag{b}{${\pi\over2}-\beta_1$}
\psfrag{A3}{$A_k$}
\psfrag{A4}{$A_{k+1}$}
\psfrag{Oi4}{$\oma{k+1}$}
\psfrag{V1}{$\tilde v_k=\ell\,\omb{k}$}
\psfrag{v1}{$\tilde v_{k}\,\cos\beta_1$}
\psfrag{V2}{$\tilde v'_k=\ell\,\omc{k}$}
\psfrag{v2}{$\tilde v'_{k}\,\cos\beta_1$}
\includegraphics[height=10cm]{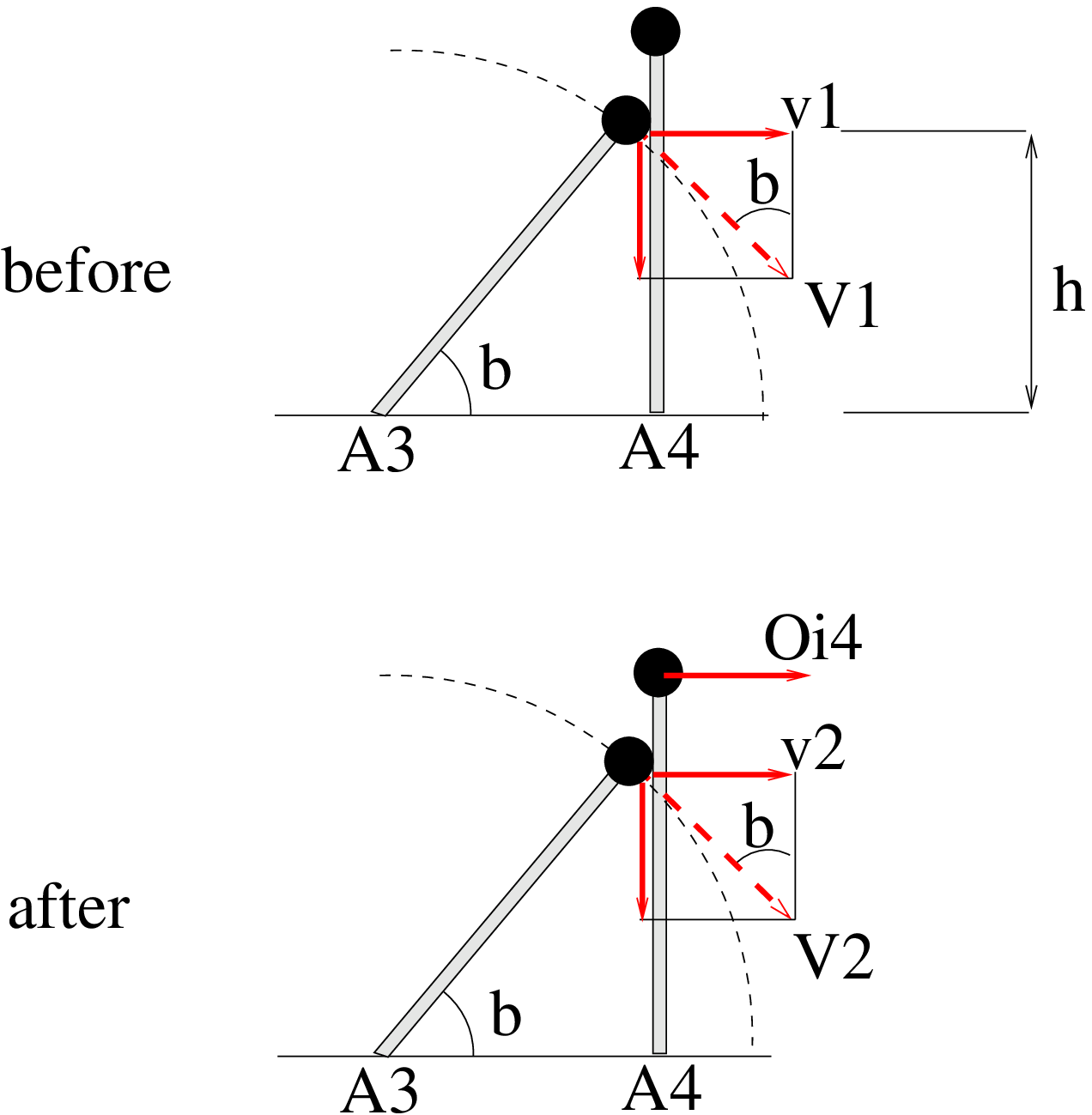}
\caption{The  collision between the \kth{k} and \kth{(k+1)} rods.}
\label{fig:rod-collision}
\end{center}
\end{figure}

Next we apply conservation of angular momentum with respect to point $A_{k+1}$.
Just before the collision the \kth{k} rod has angular velocity $\omb{k}$
and thus translational velocity $\tilde v_k=\ell \omb{k}$. From figure 4 one
can see that only the component $\tilde v_k\cos\beta_1$
contributes to the angular momentum calculated around the point $A_{k+1}$, and that it
does so with impact parameter $\ell\cos\beta_1$.
The \kth{(k+1)} rod has no angular momentum initially.
Therefore
\bb
  L_{initial}\= m\, (\tilde v_k\cos\beta_1) \, (\ell\cos\beta_1)
     \= m\,\ell^2 \, \omb{k} \, \cos^2\beta_1~.
\ee
After the collision the  \kth{(k+1)}  rod rotates around the point $A_{k+1}$ with
angular velocity $\oma{k+1}$ and thus has angular momentum
$I\oma{k+1}$.
The \kth{k} rod has the new angular velocity $\omc{k}$, again around
the point $A_k$. Therefore it will contribute an angular
momentum $m\,\ell^2 \, \omc{k} \, \cos^2\beta_1$ with respect
to $A_{k+1}$.
Therefore
\bb
  L_{final}\=m\,\ell^2\,\omc{k}\cos^2\beta_1+I\oma{k+1}~.
\ee
Conservation of angular momentum ($L_{initial}=L_{final}$) yields
\beq
  \omb{k} \, \cos^2\beta_1\= \omc{k}\cos^2\beta_1+\oma{k+1}~.
\label{eq:6}
\eeq
The system of  equations \calle{eq:7} and \calle{eq:6} can be solved easily
for $\oma{k+1}$ and $\omc{k}$:
\beqn
   \oma{k+1} &=& f_+\, \omb{k} ~,
\label{eq:2} \\
   \omc{k} &=& {\oma{k+1}\over f_-} ~, \nonumber
\eeqn
where
$$
   f_\pm\eq{2\over\cos^2\beta_1\pm 1/\cos^2\beta_1}~.
$$

Now consider the \kth{k} rod as it falls from the vertical to angle $\beta_1$, its position just
before the collision.
Conservation of total energy yields
\bb
  {1\over2}I\oma{k}^2 + mg\ell \=
  {1\over2}I\omb{k}^2 + mg\ell\,\cos\beta_1~,
\ee
or
\beq
   \omb{k}^2\= \oma{k}^2 + {2g\over\ell}\,(1-\cos\beta_1)~.
\label{eq:3}
\eeq
Combining equations \calle{eq:2} and \calle{eq:3} we find:
\beq
  \oma{k+1}^2 \= f_+^2\, \oma{k}^2 + b
\label{eq:4}
\eeq
where
$$
   b\= {2g\over\ell}\,f_+^2\,(1-\cos\beta_1)~.
$$
Equation \calle{eq:4} is a mixed progression
(i.e., a combination of an arithmetic
and a geometric progression) and can be solved by well-known techniques (see
appendix \ref{sec:mixed}). The result is
\bb
  \oma{k}^2\= f_+^{2(k-1)}\, \oma{1}^2
              + b\,{1-f_+^{2(k-1)}\over 1-f_+^2}~.
\ee
Recall that $\oma{1}$ is the initial angular velocity of the first rod
caused by the initial external push.

We show now that $f_+<1$. Since $\beta_1\ne0,\pi/2$,
$x=\cos^2
\beta_1\ne1,0$.
Then $(x-1/x)^2>0\Rightarrow x^2+1/x^2-2>0\Rightarrow x^2+1/x^2>2$. From the
last inequality it follows that $f_+=2/(x^2+1/x^2)<1$.

Since $f_+<1$ it follows that $\lim_{n\to+\infty}f_+^n=0$ and therefore
$$
 \lim_{k\to+\infty} \oma{k}^2\=  {2g\over\ell}\,(1-\cos\beta_1)\,
                                  {f_+^2\over 1-f_+^2}\eq\Omega^2~.
$$
Thus deep into the chain we find translational invariance:
the initial angular velocity imparted to a rod by its
neighbor becomes independent of position. Notice that in this limit the initial
push given to the first rod becomes irrelevant.

\subsection{Wave Speed}

We can obtain the limiting speed of the wave by computing the time
between collisions, working well into the chain where this becomes
independent of position.

Apply conservation of energy for the \kth{n} rod as it begins moving
and after it falls through an arbitrary angle $\theta$. This yields:
$$
 {1\over2}I\oma{n}^2+mg\ell\={1\over2}I\omega_{n}^2+mg\ell\,\cos\theta~.
$$
Setting $\omega_n=d\theta/dt$, we can separate $t$ from $\theta$
and solve for the time required for the rod to move from $\theta=0$ to
$\theta=\beta_1$:
$$
   \int_0^{T_n} dt \= \int_0^{\beta_1}
   {d\theta\over\sqrt{\oma{n}^2+{2g\over\ell}-{2g\over\ell}\cos\theta}}~.
$$
This integral can be expressed in terms of the complete elliptic integral
of the first kind $K(k)$ (see Appendix \ref{sec:elliptic}):
$$
   T_n\={2\over\sqrt{a_n+c}}\, \lb K\lp k_n\rp-F\left({\pi-\beta_1\over2},k_n\right)\rb~,
$$
where $a_n=\oma{n}^2+{2g\over\ell}$, $c={2g\over\ell}$ and
$k_n=\sqrt{2c\over a_n+c}$.

\begin{figure}[htb!]
\begin{center}
 \psfrag{l}{$\ell$} \psfrag{G}{$G$} \psfrag{d}{$d$}
 \includegraphics[width=8cm,angle=90]{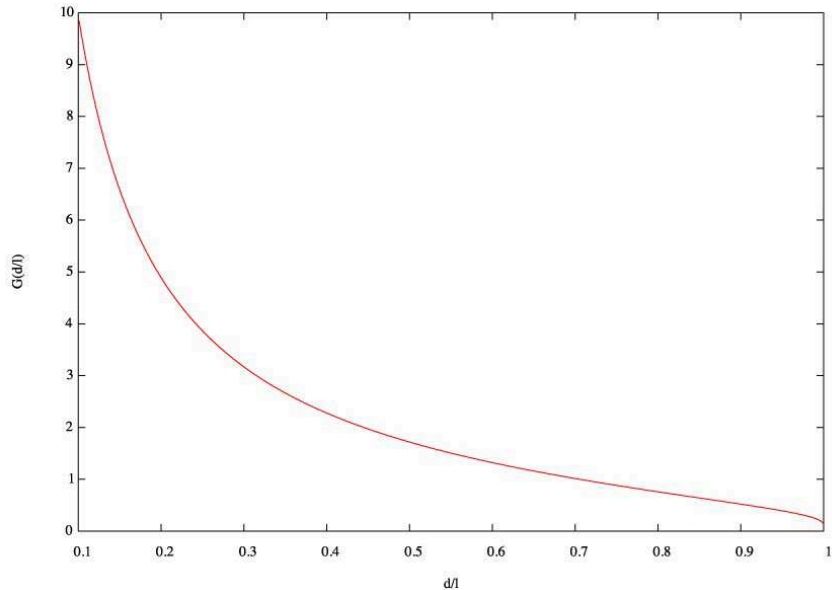}
 \caption{The scaling function $G(d/\ell)$.}
 \label{fig:graph}
\end{center}
\end{figure}

In the limit of large $n$, the time $T_n$ approaches a limiting
value
$$
   T\={2\over\sqrt{a+c}}\,
   \lb K\lp k\rp-F\left({\pi-\beta_1\over2},k\right)\rb~,
$$
where $a=\Omega^2+{2g\over\ell}$, $c={2g\over\ell}$ and
$k=\sqrt{2c\over a+c}$.
The wave therefore approaches a
limiting speed $v=d/T$ given by
$$
  v\= {d\over 2}\, {\sqrt{a+c}\over K(k)-F({\pi-\beta_1\over2},k)}~.
$$
A little algebra lets us write the wave speed in the scaling form of
equation \calle{eq:1}:
$$
    v \= \sqrt{g\ell}\,G~,
$$
with
$$
   G ({d\over\ell}) \=\frac{d}{\ell} \, \frac{1}{ k[K(k)-F({\pi-\beta_1\over2},k)]}
$$
and
$$
  k^2 \= \frac{2(1-f_+^2)}{(1-\cos\beta_1)f_+^2 +2(1-f_+^2)}~.
$$
Since $f_+$ and thus $k$ depend only on
$\beta_1=\sin^{-1}(d/\ell)$, this $G$ is indeed a function only of
$d/\ell$, as required by scaling. The scaling function $G$ is
plotted in figure \ref{fig:graph}.

No simple closed expression exists for complete elliptic integrals, but
some insight into our solution comes from looking at the limit of very closely
spaced rods ($d\ll\ell$). In this limit $\beta_1\approx d/\ell$
and $f_+^2\approx1-\beta_1^4$.
Using these in the above expression yields $k^2\approx 4\beta_1^2\ll1$.
Then
$$
 K(k)-F\left(\frac{\pi-\beta_1}{2},k\right) \=
 \int_{\frac{\pi-\beta_1}{2}}^{\frac{\pi}{2}}
 \frac{dt}{\sqrt{1-k^2 \sin^2t}}
 ~\approx~  \int_{\frac{\pi-\beta_1}{2}}^{\frac{\pi}{2}} dt
 \= \frac{\beta_1}{2}~.
$$
From this we find
$$
   G ({d\over\ell}) ~\approx~ \frac{1}{d/\ell}~.
$$
Thus the wave in very closely spaced rods moves very fast.

In a similar way one can examine the other extreme geometrical
limit, $d/\ell$ slightly smaller than unity.
Put $\beta_1=\pi/2-\sqrt{2\epsilon}$.
Then
$$
\frac{d}{\ell}=\sin\beta_1=\cos\sqrt{2\epsilon}\approx1-\epsilon
$$
while $f_+\approx 4\epsilon$ and $k^2\approx1-16\epsilon^2$.
For $k$ very near unity the complete elliptic integral is approximately
$K(k) \approx \ln(4/k')$ where $k'^2 = 1-k^2$ \cite{GR}. Here this gives
$K(k)\approx\ln(1/\epsilon)$, which diverges as $\epsilon$ approaches
zero. However in this limit $F$ is finite:
$$
  F\left(\frac{\pi-\beta_1}{2},k\right)
  = \int_0^{{\pi\over4}+{\sqrt{\epsilon}\over2}}
     \frac{dt}{\sqrt{1-k^2 \sin^2t}}
 \approx
  \int_0^{\pi\over4}
     \frac{dt}{\cos t}
 = \ln(\sqrt{2}+1)~.
$$
In the limit $d\approx\ell$ using the above approximations yields
$$
    v \approx \frac{\sqrt{g\ell}}{\ln\lb
              \frac{\ell}{(1+\sqrt{2})(\ell-d)}\rb}~,
$$
or $v \approx\sqrt{g\ell}\,G(d/\ell)$, where
 $$
    G({d\over\ell})
         ~=~ - {1\over\ln(1+\sqrt{2})+\ln(1-{d\over\ell})}~,
 $$
 as $d/\ell$ increases to 1.
  Thus as $d$ gets very close to $\ell$ the wave speed drops to
zero. Physically this occurs because one rod gives a very small
push to the next in line, which as a result takes a very long time
to fall. The reader might wish to conduct a quick experiment to
verify this conclusion or a more careful one to compare our
theoretical results with experiment. Of course, the reader has
noticed the consistency of our results with that of McLachlan et
al. \calle{eq:1}.

\section{Conclusions}

In this paper we have presented a set of assumptions for the
propagation of the domino wave and we have computed the
corresponding limiting speed. For simplicity we have presented the
solution for a simplified geometry. However, the reader can easily
transfer the solution to the case of dominoes with the shape of
parallelepipeds---with appropriate adjustments of course.

\section*{Acknowledgements}

We thank the referees for bringing to our attention the articles
of Walker \cite{Walker} and Shaw \cite{Shaw},
and the book of Banks \cite{Banks}.

\section*{Note Added in Proof}

After this work was completed a paper \cite{Leeuwen} with
somewhat similar analysis appeared on the Cornell archives.

\newpage
\appendix

\section{Mixed Progression}
\label{sec:mixed}

Consider a sequence $a_k,~k=1,2,\dots$, with the recurrence relation
$$
   a_k\= r\, a_{k-1} + b~.
$$
This is known as a mixed progression.
We want to express $a_k$ in terms of $a_1$,  $r$ and $b$.
Multiply both sides by $r^{n-k}$ and sum from $k=2$ to $n$ (with $n\ge2$):
$$
\sum_{k=2}^n r^{n-k} a_k = \sum_{k=2}^n \left( r^{n-k+1} a_{k-1}
+ r^{n-k} b \right).
$$
In the first term on the right-hand side replace $k$ by $k'=k-1$ and
in the second replace $k$ by $k'=n-k$. This yields
$$
\sum_{k=2}^n r^{n-k} a_k = \sum_{k'=1}^{n-1} r^{n-k'} a_{k'}
+ b \sum_{k'=0}^{n-2} r^{k'} .
$$
The first two sums have nearly all terms in common (all but the
\kth{n} on the left and the first on the right). Cancelling the
terms in common and evaluating the third sum yields the desired
solution:
$$
a_n = r^{n-1} a_1 + b \frac{1-r^{n-1}}{1-r}~.
$$

\section{The Elliptic Integral of First Kind}
\label{sec:elliptic}

The elliptic integral of the first kind \cite{GR,Boas} is defined
by
$$
  F(\phi_0,k)\=\int_0^{\phi_0}{d\phi\over\sqrt{1-k^2\sin^2\phi}}~,
  ~~~0\le k < 1~.
$$
When $\phi_0=\pi/2$ this is called the {complete elliptic integral
of the first kind}, denoted by $K(k)$:
$$
  K(k)\=\int_0^{\pi/2}{d\phi\over\sqrt{1-k^2\sin^2\phi}}~.
$$

The integral
$$
   I(\theta_0)\=\int_0^{\theta_0} {d\theta\over\sqrt{a-c\cos\theta}}~,
   ~~~0<c<a~,
$$
can be expressed in terms of the elliptic integral of the first kind as follows.
First make the change of variable $\theta=\pi-2t$:
$$
   I(\theta_0)\= 2\,\int_{\pi/2-\theta_0/2}^{\theta_0/2} {dt\over\sqrt{a+c\cos(2t)}}~.
$$
Using the identity $\cos{2t}=1-2\sin^2{t}$ this becomes:
\bb
   I(\theta_0)
   &=& 2\,\int_{\pi/2-\theta_0/2}^{\pi/2} {dt\over\sqrt{(a+c)-2c\sin^2t}}\\
     &=& {2\over\sqrt{a+c}}\,
      \int^{\pi/2}_{\pi/2-\theta_0/2} {dt\over\sqrt{1-{2c\over a+c}\sin^2t}}~.
\ee
Finally we set
$$
   k^2\eq{2c\over a+c}~,
$$
and we rewrite the above result in the form
\bb
  I(\theta_0)&=&{2\over\sqrt{a+c}} \left( \,\int_0^{\pi/2}{dt\over\sqrt{1-k^2\,\sin^2t}}
        -\int_0^{\pi/2-\theta_0/2}
         {dt\over\sqrt{1-k^2\,\sin^2t} }\right)\\
    &=&{2\over\sqrt{a+c}}\,\lb K(k)-F\lp{\pi-\theta_0\over2},k\rp\rb~.
\ee


\end{document}